\documentclass[12pt]{article}

\usepackage[utf8]{inputenc}
\usepackage[T1]{fontenc}
\usepackage{lmodern}
\usepackage{amsmath,amssymb,amsthm}
\usepackage[pdf,arrow,curve]{xy}
\usepackage[authoryear]{natbib}
\usepackage{hyperref}
\usepackage{geometry}
\makeindex

\theoremstyle{plain}
\newtheorem{theorem}{Theorem}[section]
\newtheorem{proposition}[theorem]{Proposition}
\newtheorem{lemma}[theorem]{Lemma}
\newtheorem{corollary}[theorem]{Corollary}

\theoremstyle{definition}
\newtheorem{definition}[theorem]{Definition}
\newtheorem{example}[theorem]{Example}

\newcommand{\imp}{\rightarrow}

\newcommand{\BImp}{\Leftrightarrow}
\newcommand{\bimp}{\leftrightarrow}

\newcommand{\proves}{\vdash}

\newcommand{\nproves}{\nvdash}

\newcommand{\lc}[1]{#1\!\!:\!\!}

\newcommand{\pa}{\mbox{\sf PA}}
\newcommand{\lp}{{\sf LP}}

\newcommand{\gn}[1]{\ulcorner {#1} \urcorner }

\title{\textbf{Provability of Consistency: A New Foundational Landscape}}
\author{Sergei Artemov\\
\small CUNY Graduate Center, New York, NY\\
\small \texttt{sartemov@gc.cuny.edu}}
\date{}

\begin{document}

\maketitle

\begin{abstract}
This paper surveys recent developments on the provability of consistency. For nearly a century, the Unprovability of Consistency Thesis (UCT) has held that sufficiently rich systems, such as Peano Arithmetic (\pa), cannot formalize a proof of their own consistency. Recent studies demonstrate that the standard G\"odelian consistency formula, $\mathsf{Con}_{\mathsf{PA}}$, on which the UCT argument is based, represents a {\bf non-equivalent redefinition} of the original Hilbertian notion of consistency. We construct a rigorous mathematical proof of the Hilbertian consistency of \pa\ that is formalizable within \pa\ itself, thereby refuting UCT. While G\"odelian consistency formulas remain indispensable for comparing the relative strengths of theories, Hilbertian consistency enables internal self-verification. Beyond vindicating Hilbert's consistency program, these findings dismantle the Lucas--Penrose argument and remove perceived logical barriers to self-verifying Automated Theorem Proving and Artificial Intelligence.

\end{abstract}

\section{Introduction}

The informal Brouwer--Heyting--Kolmogorov (BHK) semantics maps classical formulas to their constructive readings. Later foundational frameworks, such as Martin-L\"of's Intuitionistic Type Theory, bypassed the classical language altogether to construct intuitionistic analogs directly, wherein formulas and proofs gave way to dependent types and terms. This historical trajectory created the persistent illusion that the constructive domain was strictly insulated from the classical one.

Conceived by G\"odel in the 1930s and developed into the Logic of Proofs (\lp) in the 1990s, justification logic made BHK semantics mathematically rigorous, establishing an explicit bridge between the classical and constructive domains \cite{Art01a,AF19}. Within this framework, the constructive semantics of a universal statement $\forall x\, F(x)$ is realized as a \emph{selector proof}: a pair consisting of an explicit computable function $s(x)$ (the \emph{selector}) and a \emph{verifier} $v$ demonstrating that, for every $x$, $s(x)$ is a proof of $F(x)$. 

Crucially, this constructive apparatus also decisively contributes to classical foundations. Selector proofs are routinely accepted in everyday mathematical practice and machine reasoning, are conservative over Peano Arithmetic (\pa), and can be safely employed as a sound \pa\ tool. Utilizing this, we construct a mathematical proof of \pa-consistency that is formalizable in \pa, thereby refuting the long-standing paradigm that consistency is inherently unprovable within a theory itself.

\subsection{The Unprovability of Consistency Thesis}

In the 1930s, G\"odel's Second Incompleteness Theorem (G2) was rapidly interpreted as establishing the \emph{Unprovability of Consistency Thesis} (UCT), summarized in standard reference works as follows:
\begin{quote}
``There exists no consistency proof of a system that can be formalized in the system itself'' (\emph{Encyclop\ae dia Britannica}).
\end{quote}
This thesis has exerted a profound influence on the foundations of mathematics, the theory and practice of formal verification, and the philosophy of mind and cognition. 

In strict mathematical reality, G\"odel demonstrated that a \emph{specific consistency formula} ($\mathsf{Con}_{\mathsf{PA}}$) is unprovable in \pa, providing a potent yardstick for measuring the relative proof-theoretic strength of formal theories. 

Both Hilbert and G\"odel issued explicit cautions against overinterpreting G2. G\"odel voiced these reservations in his 1931 paper \cite{God31} (immediately following Theorem XI), in his 1930--1931 correspondence with von Neumann, and even decades later in private conversations with Gerald Sacks in Princeton during the early 1970s, maintaining that Hilbert's program had not been erased by G2 but remained ``very much alive and even more interesting than it initially was.'' Despite these warnings, John von Neumann drew the broader philosophical conclusion and actively promoted UCT, the view that mathematical consistency itself is unprovable internally.

We note that the standard consistency formula $\mathsf{Con}_{\mathsf{PA}}$ is a redefinition of consistency that is \emph{not equivalent in \pa} to the original Hilbertian notion of consistency \cite{Art25b}. By offering a mathematical proof of the original Hilbert consistency of \pa\ formalizable in \pa\ \cite{Art25a}, we render UCT in its naive form false. To be entirely clear: we are not correcting G\"odel's mathematics. G2 remains a flawless logical truth. Rather, we are correcting the philosophical orthodoxy promoted by von Neumann that overextended a technical theorem into an absolute barrier against self-verification.

\subsection{A Historical Precedent: Von Neumann's Quantum Impossibility Proof}

The elevation of a valid technical result into a restrictive philosophical dogma via hidden assumptions has a direct historical parallel in twentieth-century physics. As Dennis Dieks recounts regarding von Neumann's 1932 theorem against hidden variables in quantum mechanics:
\begin{quote}
``According to what has become a standard history of quantum mechanics, in 1932 von Neumann persuaded the physics community that hidden variables are impossible \dots This state of affairs lasted \dots until Bell in 1966 exposed von Neumann's proof as obviously wrong. The realization that von Neumann's proof was fallacious then rehabilitated hidden variables and made serious foundational research possible again.'' \citep{Dieks2017}
\end{quote}
Von Neumann's justification of UCT shares an identical logical structure with his quantum impossibility proof: both were mathematically sound deductions built upon overly restrictive, fallacious assumptions about physical or metamathematical concepts. Just as Bell (and Grete Hermann earlier) rehabilitated hidden variables by exposing von Neumann's unjustified linearity assumption, we rehabilitate internal consistency proofs by exposing the domain violation inherent in $\mathsf{Con}_{\mathsf{PA}}$.

\section{Hilbert's Notion of Consistency}

\subsection{Peano Arithmetic and Non-Standard Models}

Peano Arithmetic (\pa) is a formal first-order theory containing the constant $0$, the successor function $'$, and all primitive recursive functions together with their defining identities. In addition, \pa\ includes the Induction Schema, consisting of an infinite series of formulas:
\[
[\varphi(0) \wedge \forall x\,(\varphi(x) \imp \varphi(x'))] \imp \forall x\,\varphi(x)
\]
for each formula $\varphi(x)$ in the language of arithmetic. The Induction Schema serves as our primary archetype of a \emph{serial property}: an infinite, primitive recursive sequence of formulas. While \pa\ represents all conventional finitary computations, it employs no higher-order or set-theoretic principles.

The standard model of \pa\ is the set $\mathbb{N} = \{0, 1, 2, \dots\}$ equipped with standard addition and multiplication:
\begin{figure}[!h]
\begin{center}
\mbox{
\begin{xy}
(-60,-4)*{0};
%(-63,0)*{\mbox{$p$}};
(-60,0)*+{\bullet}="2"; 
(-50,-4)*{1};
(-50,0)*+{\bullet}="3"; 
(-40,-4)*{2};
(-40,0)*+{\bullet}="4"; 
(-30,0)*+{\ldots}="5"; 
{\ar "2";"3"};
{\ar "3";"4"};
{\ar "4";"5"};
\end{xy}
.}\end{center}
%\caption{Model $\mathcal{M}_1$.}
\end{figure}

By the Compactness Theorem, there exist countable models of \pa\ not isomorphic to the standard model. Introducing a fresh individual constant $c$, consider the theory
\[
\widetilde{\pa} = \pa \cup \{ c > 0,\ c > 1,\ c > 2,\ \dots \}.
\]
Every finite subsystem of $\widetilde{\pa}$ is consistent, so $\widetilde{\pa}$ possesses a model $\mathcal{M}$ containing ``infinite'' (non-standard) elements alongside the standard copy of $\mathbb{N}$:
\[
\begin{xy}
(-60,-4)*{0};
(-60,0)*+{\bullet}="2"; 
(-50,-4)*{1};
(-50,0)*+{\bullet}="3"; 
(-40,-4)*{2};
(-40,0)*+{\bullet}="4"; 
(-30,0)*+{\ldots}="5"; 
(-10,-4)*{c};
(-10,0)*+{\bullet}="c"; 
(-20,0)*+{\ldots}="6"; 
(0,0)*+{\ldots}="7"; 
{\ar "2";"3"};
{\ar "3";"4"};
{\ar "4";"5"};
{\ar "6";"c"};
{\ar "c";"7"};
\end{xy}
.\]

\subsection{Numerals vs. Natural Numbers}

For the purposes of Hilbert's consistency program, standard natural numbers are represented constructively as \pa-numerals (analogous to Hilbert strokes, or Zermelo, von Neumann, and Church numerals):
\[
0,\ 0^{\prime },\ 0^{\prime \prime },\ 0^{\prime \prime \prime },\ \dots .\]

Crucially, the set of standard natural numbers is not definable in \pa. G\"odel numbering assigns numerical codes to 
syntactic objects---terms, formulas, and finite sequences of formulas---using numerals.

The principal divergence between informal contentual arithmetic and formal \pa\ lies in quantification. Because standard natural numbers cannot be isolated by a predicate in \pa, the informal quantifier 
\[\mbox{\emph{``for all natural numbers $n\dots$''}}\] is translated in \pa\ into the formal first-order quantifier $\forall x\dots$, which ranges over all elements of a model, including non-standard elements. Consequently, $\forall x\,\varphi(x)$ asserts a stronger claim than 
``for all natural numbers $n$, $\varphi(n)$ holds.'' 
%\textbf{G\"odel's Second Incompleteness Theorem lives precisely within this gap.}

\subsection{Hilbert's Formulation of Consistency}
For Hilbert, the domain of finitary number theory consists of numerals. As Hilbert stated in 1928, a finitary general proposition is ``a hypothetical judgment that comes to assert something when \emph{a numeral} is given'' \citep{Hilbert1928}.
The mathematical statement of consistency for a formal theory,
\begin{equation}\label{eq:hilbert_orig}
\mbox{\it no finite sequence $S$ of formulas is a derivation of a contradiction,}
\end{equation}
becomes, upon standard G\"odel coding of finite sequences by numerals, a Hilbertian finitary general proposition:
\begin{equation}\label{eq:hilbert_numeral}
\mbox{\it any numeral $n$ is not a code of a derivation of a contradiction.}
\end{equation}
Thus, Hilbert's definition of consistency concerns exclusively finite derivations represented by their numeral codes. Non-standard elements in non-standard models of \pa---none of which can ever encode a finite sequence of formulas---are entirely irrelevant to the consistency of the theory. This observation will be substantiated by solid mathematical arguments below. 

\section{Formalization and the Strengthening Fallacy}

\subsection{Modes of Formalization and Verification}\label{Stformalization}
We distinguish two levels of formalizing mathematical statements in \pa:
\begin{enumerate}
\item \textbf{Direct formalization:} Given an informal arithmetical property $H$, a formalization procedure $f$ preserves primitive recursive operations and logical connectives while replacing ``for any natural number $n$'' with $\forall x$, yielding a \pa-formula (or set of formulas) $f(H)$.
\item \textbf{G\"odelian formalization:} Standard arithmetization assigns numerical codes to syntactic objects and converts syntactic operations into primitive recursive functions on codes, followed by direct formalization in \pa. Reasoning about mathematical objects is thus converted step-by-step into formal derivations in \pa\ concerning G\"odel codes. All formalizations in our consistency proof are standard G\"odelian.
\end{enumerate}
When investigating whether a contentual property $H$ of natural numbers holds via its formalized counterpart $f(H)$, formalization serves two roles:
\begin{itemize}
\item \textbf{Verification:} If $\pa \proves f(H)$, then $H$ holds; if $\pa \proves \neg f(H)$, then $H$ fails. If neither is provable, single-formula formalization leaves $H$ undecided.
\item \textbf{Proof and Assumption Checking:} If a contentual mathematical proof is formalized in \pa, this certifies that the proof relies on no principles outside \pa.
\end{itemize}
\subsection{Avoiding the `Vacuous Consistency' Trap}
A common objection asks: 
\begin{quote}\emph{What is the point of establishing the provability of consistency of a theory $T$ within $T$ itself, given that if $T$ were inconsistent, any statement would be vacuously provable in $T$?}
\end{quote}
This objection rests on a misunderstanding of Hilbert's program. 

\[ \mbox{``Establishing the consistency of $T$ within $T$''} \]
is a brief version of 
\[ \mbox{``Establishing the consistency of $T$ by means formalizable in $T$,''}\]
which consists of two essential components:
\begin{enumerate}
\item[(i)] A contentual mathematical proof of $T$-consistency (which fails if $T$ is inconsistent);
\item[(ii)] A G\"odelian formalization of (i) inside $T$ for assumption checking, certifying that no external axioms were smuggled into (i).
\end{enumerate}
The contentual component (i) is indispensable for \emph{any} presentation of consistency. Furthermore, our approach adds a crucial third benefit:
\begin{enumerate}
\item[(iii)] An explicit, finitely checkable proof object $\mathcal{C}$ of consistency. Previously, trusting the consistency of $T$ required trusting an infinite set of $T$-derivations; with $\mathcal{C}$, trust reduces to verifying a single finite certificate, manually or via an independent proof checker.
\end{enumerate}

\subsection{The Domain Violation and the Strengthening Fallacy}\label{fallacy}
%The semi-formal analysis of UCT and its justifications in this section serves as an explanatory supplement for the reader and is not essential to our main results on consistency proofs. 

Let $T$ be a theory. Let also $$\lc{x}_{\small T}y$$ denote the standard primitive recursive proof predicate $\mathsf{Prf}_{\small T}(x,y)$ (``$x$ is the G\"odel number of a $T$-derivation of the formula with G\"odel number $y$''). We omit `$T$' when $T$ is \pa.

We adopt the usual simplifying notational convention that allows us not to distinguish within ``:'' numerals, standard numbers, formulas, and their G\"odel numbers in a way that respects the usual understanding of free variables, when safe. In particular, we will freely use notations $\lc{n}m$, $\lc{t}F$, $\lc{t(x)}F(x)$, etc. 

Let $\mathit{Consistent}_{\small T}(x)$ stand for $\neg\lc{x}_T(0=1)$. 
Arithmetized \pa-consistency is the contentual property:
\begin{equation}\label{eq:acon}
\mbox{\it for any numeral $n$, $\mathit{Consistent}(n)$.}
\end{equation}
Notice that \eqref{eq:acon} is not a sentence in the language of \pa\ because quantification restricted to standard numerals is not definable in \pa.
The traditional approach ignores this domain restriction and defines consistency inside \pa\ as the $\Pi_1$-sentence:
\begin{equation}\label{eq:conpa}
\mathsf{Con}_{\mathsf{PA}} \equiv \forall x\, \mathit{Consistent}(x) \equiv \forall x[\neg\lc{x}(0=1)].
\end{equation}
By G\"odel's Second Incompleteness Theorem, if \pa\ is consistent, then $\pa \nproves \mathsf{Con}_{\mathsf{PA}}$. Consequently, there exists a non-standard model $\mathcal{M}$ of \pa\ in which $\mathsf{Con}_{\mathsf{PA}}$ is false. In $\mathcal{M}$, there exists an element $c$ satisfying $\lc{c}(0=1)$---G\"odel's ``monster'' proof of a contradiction. However, $c$ cannot be any standard numeral $n = 0, 1, 2, \dots$; it is a non-standard number that does not correspond to any real finite \pa-derivation.
Conceptually, $\mathsf{Con}_{\mathsf{PA}}$ is infested with non-standard ``junk proof codes'' that \pa\ is too weak to filter out. Mathematically, this yields the following proposition \cite{Art25b}:
\begin{proposition}\label{thm:strengthening}
\pa\ cannot prove the equivalence of $\mathsf{Con}_{\mathsf{PA}}$ and the contentual consistency property \eqref{eq:acon}. Specifically, $\mathsf{Con}_{\mathsf{PA}}$ is strictly stronger over \pa\ than \pa-consistency.
\end{proposition}
The proof is a straightforward combinatorial argument immediately formalizable in \pa.

\subsection{What is wrong with the standard justification of UCT?}
The semi-formal analysis in this subsection serves as an explanatory supplement for the reader and is not essential to our main results on consistency proofs. 

Deriving the Unprovability of Consistency Thesis from the unprovability of $\mathsf{Con}_{\mathsf{PA}}$ commits a \emph{strengthening fallacy}. The standard syllogism for UCT runs:
\begin{enumerate}
\item $\pa \nproves \mathsf{Con}_{\mathsf{PA}}$ (by G2);
\item \pa-consistency implies $\mathsf{Con}_{\mathsf{PA}}$;
\item Therefore, \pa\ does not prove \pa-consistency.
\end{enumerate}
For (3) to follow within \pa, step (2) would need to be provable in \pa, which fails by Proposition~\ref{thm:strengthening}. 
A more sophisticated defense of UCT replaces (2) with the \emph{Formalization Principle}:
\begin{enumerate}
\item[$2'$.] A contentual arithmetical proof of a property $H$ formalizes as a \pa-derivation of the standard formalization $f(H)$ (cf. Section~\ref{Stformalization}).
\end{enumerate}
Since $\mathsf{Con}_{\mathsf{PA}}$ is the standard formalization of the property ``\pa\ is consistent,'' ($2^\prime$) yields that arithmetic does not prove \pa-consistency.

In \cite{Art25a}, we refute ($2'$) directly by exhibiting a contentual Hilbert-style proof of \pa-consistency, formalizable in \pa, that does not yield a \pa-proof of $\mathsf{Con}_{\mathsf{PA}}$. To be precise, our proof of \pa-consistency is formalizable in \pa, but not as a derivation of $\mathsf{Con}_{\mathsf{PA}}$. Intrinsically, our proof is a generalized (and successful) version of Hilbert's $\epsilon$-substitution method {\bf which did not aim at proving $\mathsf{Con}_{\mathsf{PA}}$}. For a more in-depth analysis, cf. Section~\ref{noncompact}.

\section{Serial Properties and Selector Proofs}

\subsection{Serial Properties}
\begin{definition}
A \emph{serial property} is a primitive recursive sequence of formulas $F_0, F_1, F_2, \dots$, {\bf finitely presented} by a primitive recursive enumerating \pa-term.
\end{definition}
Serial properties are ubiquitous across logic and mathematics: the set of terms, formulas, and tautologies of \pa; the Induction Schema of \pa; the axiom set of \pa\ itself; reflection schemas; and comprehension/separation schemas in $\mathsf{ZF}$. Each is finitely represented by a primitive recursive term.
\begin{definition}
The \emph{consistency schema} $\mathsf{ConS}_{\mathsf{PA}}$ is the serial property:
\begin{equation}\label{eq:cons_schema}
\mathit{Consistent}(n), \quad n = 0, 1, 2, \dots
\end{equation}
\end{definition}
\begin{proposition}\label{prop:equiv}
The consistency schema $\mathsf{ConS}_{\mathsf{PA}}$ is provably equivalent in \pa\ to Hilbertian consistency \eqref{eq:acon}.
\end{proposition}
\begin{proof}
This follows by straightforward combinatorial reasoning regarding G\"odel numbering of finite sequences of formulas, which is naturally G\"odel-formalizable in \pa.
\end{proof}
Unlike $\mathsf{Con}_{\mathsf{PA}}$, which is a $\Pi_1$-formula distorted by an unbounded universal quantifier, $\mathsf{ConS}_{\mathsf{PA}}$ is presented by a quantifier-free primitive recursive term and provides a faithful {\bf finite}  representation of consistency in \pa.

The consistency schema $\mathsf{ConS}_{\mathsf{PA}}$ is known to be instance-provable in \pa: 
\[ \pa\proves\mathit{Consistent}(n),\ \mbox{\emph{ for each $n = 0, 1, 2, \dots$},} \]
\[  \pa+\mathsf{ConS}_{\mathsf{PA}} = \pa,\]
hence $\pa+\mathsf{ConS}_{\mathsf{PA}}$ is strictly weaker than $\pa+\mathsf{Con}_{\mathsf{PA}}$. We refer to this fact as $\mathsf{ConS}_{\mathsf{PA}}$ is strictly weaker than $\mathsf{Con}_{\mathsf{PA}}$ in \pa, or 
\[ \mbox{\emph{$\mathsf{Con}_{\mathsf{PA}}$ is strictly stronger than \pa-consistency over \pa. }}\] 

\subsection{Selector Proofs and Conservativity}
Proofs of serial properties take the form of \emph{selector proofs}, a concept traceable to Hilbert's $\epsilon$-substitution method, Brouwer, Heyting, Kolmogorov, and Kreisel, and also appear as the result of the Artemov--Yavorskaya realization of the first-order intuitionistic logic into the First-Order Logic of Proofs ($\mathsf{FOLP}$) \citep{AF19}. 
Specifically, a selector proof emerges naturally through a canonical two-step logical transformation (cf. \cite{AF19}, Chapter 10):

Starting from a constructive universal statement $\forall x\, F(x)$, the G\"odel translation embeds it into First-Order $\mathsf{S4}$ ($\mathsf{FOS4}$) as $\Box \forall x \Box F(x)$, asserting that it is provable that for every $x$, $F(x)$ is provable. Applying the Artemov--Yavorskaya realization theorem then converts each provability modality $\Box$ in $\mathsf{FOS4}$ into an explicit proof term: the inner modality $\Box F(x)$ is realized by a computable selector term $s(x)$ parameterized by $x$, while the outer modality is realized by a finite verifier derivation $v$, producing the exact pair $\langle s(x), v \rangle$: 
$$\lc{v} \forall x [\lc{s(x)}F(x)]. $$

\begin{definition} Let $\mathcal{F}$ be a serial property ${F_0, F_1, \dots, F_n, \dots}$, and $F(x)$ be a natural arithmetical term for a primitive recursive function that, given $n$, computes the G\"odel number of $F_n$. 
A \emph{selector proof} of a serial property $\mathcal{F} = {F_0, F_1, \dots, F_n, \dots}$ in a theory $T$ is a pair $\langle s(x), v \rangle$ consisting of:
\begin{enumerate}
\item[(i)] A \emph{selector} $s(x)$: a primitive recursive term in $T$ that, given a numeral $n$, computes the G\"odel number of a $T$-derivation of $F_n$;
\item[(ii)] A \emph{verifier} $v$: a conventional finite derivation in $T$ establishing that the selector succeeds for all $x$:

\[T \proves \forall x\,[\lc{s(x)} F(x)].\]

\end{enumerate}
\end{definition}
Note that Hilbert's guidelines for consistency proofs \citep{Zach2007} explicitly anticipate selector proofs: 
\begin{quote}
\emph{Hilbert called for an effective operation transforming any formal derivation into a special form, accompanied by proofs that the operation always succeeds and that derivations of the special form cannot prove a contradiction.}
\end{quote}
\begin{proposition}\label{prop:conservativity}
If a serial property $\mathcal{F} = {F_0, F_1, \dots, F_n, \dots}$ is selector-provable in \pa, then each instance $F_n$ is conventionally provable in \pa.
\end{proposition}
\begin{proof}
We invoke the \emph{explicit reflection principle}, $$\lc{t}F \imp F,$$ whose provability in \pa\ was suggested by G\"odel in 1938 \cite{God38} and independently proven by Artemov and Strassen in 1992 \cite{AS92}, cf. also \cite{Art01a,AF19}. Suppose $$\pa \proves \forall x[\lc{s(x)} F(x)].$$ Instantiating with any standard numeral $n$ yields $$\pa \proves \lc{s(n)} F(n).$$ Applying explicit reflection in \pa, we obtain $$\pa \proves F_n.$$
\end{proof}
Thus, selector proofs add no new theorems over \pa:  $$\pa + \text{selector proofs} = \pa;$$ they are simply finite \pa-derivations that uniformly certify infinite families of \pa-theorems.
\begin{example}[Complete Induction]
The Complete Induction principle ($\mathcal{CI}$) states that for any formula $\psi$, \[\mbox{if $\forall x[\forall y < x,\psi(y) \imp \psi(x)]$, then $\forall x\psi(x)$.}\] The standard textbook proof applies ordinary \pa-induction to the auxiliary formula $\forall y < x,\psi(y)$ to derive $\mathcal{CI}(\psi)$. Formalized, this gives a primitive recursive selector $s(\gn{\psi})$ that builds a \pa-derivation of $\mathcal{CI}(\psi)$ for every $\psi$, verified in \pa.
\end{example}
\begin{example}[Multiplying polynomials] {\em The product of polynomials is a polynomial.}\footnote{A polynomial is a term $a_nx^n+a_{n-1}x^{n-1}+\ldots+a_1x+a_0$, where each $a_i$ is a numeral and $x$ a variable. In $f\cdot g$, ``$\cdot$'' stands for the usual \pa\ multiplication.}

Here is its standard mathematical proof.
Given a pair of polynomials $f,g$, using the well-known formula, calculate coefficients of the product polynomial $p_{f\cdot g}$, and {\bf prove in arithmetic} that
\begin{equation}\label{f.g}
f\!\cdot\! g = p_{f\cdot g}.
\end{equation}
This is a selector proof: for each $f,g$, it finds a proof of (\ref{f.g}) in \pa.
\end{example}
\begin{example}[Propositional Tautologies]
For any formula $X$, the double negation law $$X \bimp \neg\neg X$$ is proved in \pa\ by a selector $D(\gn{X})$ that constructs the standard propositional derivation tailored to the syntactic structure of $X$, accompanied by a \pa-verifier confirming that $D(\gn{X})$ is always a valid derivation. The same argument works for any tautology containing a propositional variable. 
\end{example}
These examples illustrate that most of metamathematics, which is assumed to be formalizable in \pa, uses selector proofs, since it operates on finite combinatorial objects rather than on integers. 
\begin{example}[Proving Self-Consistency of $\pa^\omega$]
Define \[\pa_0 = \pa,\ \ \pa_{i+1} = \pa_i + \mathsf{Con}_{\mathsf{PA}_i},\ \ \mbox{and}\ \  \pa^\omega = \bigcup_{i=0}^\infty \pa_i.\] Given any derivation $D$ in $\pa^\omega$, one can effectively find $n$ such that $D$ is a derivation in $\pa_n$; since $\mathsf{Con}_{\mathsf{PA}_n}$ is an axiom of $\pa^\omega$, $D$ cannot end with $0=1$. This is a valid selector proof of the consistency of $\pa^\omega$ formalizable in $\pa^\omega$, even though the G\"odelian consistency remains unprovable: $\pa^\omega \nproves \mathsf{Con}_{\mathsf{PA}^\omega}$ by G2.
\end{example}
How has mathematical logic lived without formally acknowledging selector proofs? Just as mathematics in general lived before the 1930s without formally acknowledging computability and instead operated with specific algorithms (e.g., the Euclidean Algorithm) on a case-by-case basis. At some point, the need to study computability/noncomputability in general emerged, which cannot be done without a rigorous mathematical definition of computability. Likewise, in foundations, logic lives using the instance provability of serial properties of interest, with the intuitive understanding that these proofs can be supported by verification with given means. However, since the need to study provability/unprovability in a theory has emerged within the Hilbert consistency program, we must develop a general, rigorous theory of selector proofs, as they are ubiquitous in mathematics, metamathematics, and machine reasoning. 

The analysis of foundational boundaries was incomplete because logic lacked a rigorous formal object for uniform instance-provability. Formalized selector proofs let us separate non-compact serial provability from compact object-level provability and offer a deeper analysis of provability.  

\section{Main Result: A Proof of \pa-Consistency in \pa}
\subsection{Partial Truth Definitions}
To prove that no \pa-derivation $S$ ends with $0\!=\!1$, the natural semantic argument asserts that every formula in $S$ is true in the standard model while $0\!=\!1$ is false. By Tarski's Undefinability Theorem, a single universal truth predicate for all arithmetical formulas is not definable in \pa. This makes the standard model consistency proof not formalizable in \pa.

However, proof theory provides \emph{partial truth predicates} $\mathit{Tr}_n(x, y)$ defined primitive recursively for $\Sigma_n$-formulas of bounded quantifier depth $n$ (cf. \cite{Art25a}), which help to carry out the above argument on the derivation-by-derivation basis with the subsequent uniform verification of this procedure as a finite proof in \pa. 
\begin{proposition}[Properties of Partial Truth Definitions]\label{prop:partial_truth}
For each $n = 0, 1, 2, \dots$, there exists a primitive recursive construction of a $\Sigma_{n}$-formula $\mathit{Tr}_n(\gn{\varphi}, y)$ (where $y$ encodes parameter assignments) satisfying:
\begin{enumerate}
\item $\mathit{Tr}_n(\gn{\varphi}, y)$ obeys Tarski's compositional truth clauses for Boolean connectives, bounded and unbounded quantifiers, and \emph{Modus Ponens} for all $\varphi \in \Sigma_n$.
\item For every $\Sigma_n$-formula $\varphi$, \pa\ proves Tarski's biconditional:
\[\mathit{Tr}_n(\gn{\varphi}, y) \BImp \varphi(y).\]
In particular, $\neg\mathit{Tr}_n(\gn{0\!=\!1}, y)$ is provable in \pa.
\item $\mathit{Tr}_n(\gn{A}, y)$ is provable in \pa\ for every axiom $A$ of \pa\ of syntactic depth $\leq n$.
\end{enumerate}
All proofs of these properties are rigorous contentual arguments requiring no metamathematical consistency assumptions about \pa.
\end{proposition}
\subsection{Construction of the Consistency Proof}
\begin{theorem}[Main Theorem]\label{thm:main}
The \pa-consistency represented by the consistency schema $\mathsf{ConS}_{\mathsf{PA}}$ admits a contentual mathematical proof that is G\"odel-formalizable as a selector proof in \pa:
\[\pa \proves \forall x\,[\lc{s(x)} \mathit{Consistent}(x)]\]
for an explicit primitive recursive selector term $s(x)$.
\end{theorem}
\begin{proof}[Proof Sketch]
Given any numeral $d$ coding a finite sequence $S$ of formulas that constitutes a valid \pa-derivation, we primitive recursively compute an upper bound $n = n(d)$ on the syntactic complexity ($\Sigma_n$-depth) of all formulas occurring in $S$.
We then argue by induction on the length of the derivation $S$ that every formula $\varphi$ in $S$ (under any parameter evaluation $y$) satisfies $\mathit{Tr}_n(\gn{\varphi}, y)$. By Proposition~\ref{prop:partial_truth}, all axioms of \pa\ occurring in $S$ satisfy $\mathit{Tr}_n$, and all rules of inference of \pa\ preserve $\mathit{Tr}_n$. Consequently, $\mathit{Tr}_n$ acts as a syntactic invariant across the entire derivation $S$. Because $\neg\mathit{Tr}_n(\gn{0\!=\!1}, y)$ holds, the formula $0\!=\!1$ cannot appear in $S$, establishing $\mathit{Consistent}(d)$.

Since all steps---computing $n(d)$, constructing the partial truth definition $\mathit{Tr}_n$, verifying axiom validity and rule preservation, and deducing $\neg\mathit{Tr}_n(\gn{0=1}, y)$---are primitive recursive, this contentual proof translates directly into a primitive recursive selector term $s(x)$ and a finite verifier derivation in \pa\ proving $\forall x[\lc{s(x)}\mathit{Consistent}(x)]$.
\end{proof}

Theorem ~\ref{thm:main} complements the Mostowski Reflexivity Theorem (MRT), which states that \pa\ proves the consistency of each of its finitely axiomatized fragments. Whereas MRT relaxes consistency into an implicit provability claim (requiring an unprovable reflection principle to detach), our selector proof retains explicit proof terms and detaches via Explicit Reflection and Proposition~\ref{prop:conservativity}.

\subsection{Contentual Filters Against Gerrymandered Selectors}
Just as traditional G2 literature uses contentual analysis to exclude contrived consistency formulas (such as Rosser's consistency predicate, which is provable in \pa\ via a syntactic trick), our two-stage methodology requires both a sound contentual proof (i) and its \pa-formalization (ii) to rule out ``gerrymandered'' selectors.

Consider, for instance, a pathological selector $p(x)$ defined as follows: given $x$, test whether $\lc{x}\!(0\!=\!1)$. If yes, let $p(x)$ be a trivial derivation of $\neg\lc{x}\!(0\!=\!1)$ from the contradiction $0\!=\!1$; if no, let $p(x)$ be the standard $\Sigma_1$-completeness derivation of $\neg\lc{x}(0\!=\!1)$. While $\pa \proves \forall x[\lc{p(x)}\neg\lc{x}(0\!=\!1)]$, $p(x)$ does not correspond to a genuine contentual consistency proof but rather demonstrates, similar to the proof of Rosser's consistency, how to get the formal answer even when $x$ is inconsistent. In contrast, our partial truth selector $s(x)$ is extracted directly from a sound, non-circular contentual proof (Theorem~\ref{thm:main}).

One might object that requiring a contentual filter introduces an informal, metamathematical judgment into what purports to be an internal consistency proof. However, this intensional requirement is not unique to selector proofs; it is equally indispensable to the traditional G\"odelian framework. In standard metamathematics, asserting that $\mathrm{PA} \nvdash \mathrm{Con}_{\mathrm{PA}}$ establishes the unprovability of consistency likewise relies on an informal contentual judgment that $\mathsf{Con}_{\mathsf{PA}}$---rather than a contrived Rosser-style predicate---genuinely expresses consistency. In both frameworks, formal syntax alone cannot distinguish natural mathematical concepts from pathological encodings. Our two-stage methodology simply makes this standard intensional discipline explicit: the formal verifier \(v\) certifies that no hidden axioms are used, while the contentual proof (i) guarantees that the selector $s(x)$ reflects genuine semantic invariance rather than a syntactic trick.

\section{Compact vs. Non-Compact Proofs}\label{noncompact}
\begin{definition}
A selector proof $\pa \proves \forall x [\lc{s(x)} F(x)]$ is \emph{compact} if there exists a finitely axiomatized subsystem $W \subset \pa$ such that all derivations $s(n)$ take place within $W$:
\[ \pa \proves \forall x\,[\lc{s(x)}_W F(x)].\]
Otherwise, the selector proof is called \emph{non-compact} \citep{Art25c, Gadsby2025}.
\end{definition}
\begin{lemma}[Gadsby Lemma \citep{Gadsby2025}]
Because \pa\ proves the uniform reflexivity of each of its finitely axiomatized fragments $W$, any compact selector proof $\pa \proves \forall x[\lc{s(x)}_W F(x)]$ reduces to conventional universal provability:
\[ \pa \proves \forall x F(x).\]
\end{lemma}
This lemma answers Lev Beklemishev's 2019 question\footnote{In a private communication to the author.} regarding precisely which consistency proofs are prohibited by G\"odel's Second Incompleteness Theorem. 
\begin{corollary} G2 bans compact proofs of \pa-consistency in \pa.
\end{corollary}
\par\noindent
Indeed, a compact selector proof of $\mathsf{ConS}_{\mathsf{PA}}$ would imply $\pa \proves \forall x\ \mathit{Consistent}(x)$, i.e., $\pa \proves \mathsf{Con}_{\mathsf{PA}}$. Thus, G2 does not establish the unprovability of consistency; rather, it reveals that \pa-consistency is a non-compact serial property that cannot be compressed into a single finitely provable formula.

Since proving $\mathit{Consistent}(d)$ for derivations $d$ of unbounded formula complexity requires $\Sigma_{n}$-induction for arbitrarily large $n$, our selector proof of $\mathsf{ConS}_{\mathsf{PA}}$ (like the selector proofs of Complete Induction and $\pa^\omega$-consistency) is inherently \emph{non-compact}. 

Regarding how far internal consistency proofs can be pushed, Kurahashi and Sinclaire (2019), in independent private communications to the author, noticed that \pa\ cannot selector-prove the consistency of any theory $T \supseteq \pa + \mathsf{Con}_{\mathsf{PA}}$. Technically, this observation also immediately follows from some of the results in \cite{Ign94}. However, Freund, Pakhomov, and Gadsby \cite{FP20,Gadsby2025} demonstrated that \pa\ \emph{can} selector-prove the consistency of certain proper extensions of \pa, such as $$T = \pa + \text{slow consistency of }\pa.$$

\section{Foundational and Practical Implications}

\subsection{A Complementary Dual-Aspect Framework}
Rather than forcing a choice between G\"odel and Hilbert, our results establish a balanced, complementary foundation:
%{\small
%\begin{table}[h!]
%\centering
%\renewcommand{\arraystretch}{1.3}
%\begin{tabular}{|l|p{5.2cm}|p{5.2cm}|}
%\hline
%\textbf{Attribute} & \textbf{G"odelian $\mathsf{Con}{\mathsf{PA}}$} & \textbf{Hilbertian $\mathsf{ConS}{\mathsf{PA}}$} \
%\hline
%\textbf{Finite Representation} & Single $\Pi_1$-formula containing an unbounded quantifier $\forall x$ & Quantifier-free primitive recursive enumerating term (schema) \
%\hline
%\textbf{Equivalent to Consistency in \pa?} & \textbf{No} (strictly stronger due to non-standard models) & \textbf{Yes} (provably equivalent in \pa) \
%\hline
%\textbf{Provable in \pa?} & \textbf{No} (blocked by the G2 barrier) & \textbf{Yes} (via a non-compact selector proof) \
%\hline
%\textbf{Primary Utility} & \emph{Metamathematical telescope}: comparing relative strength of theories & \emph{Syntactic microscope}: internal self-verification \
%\hline
%\end{tabular}
%\caption{Comparison of G"odelian and Hilbertian formalizations of consistency.}
%\label{tab:comparison}
%\end{table}
%}
G\"odelian consistency ($\mathsf{Con}_{\mathsf{PA}}$) operates as a \emph{metamathematical telescope}, measuring and comparing the relative proof-theoretic strengths of formal systems. Hilbertian serial consistency ($\mathsf{ConS}_{\mathsf{PA}}$) operates as a \emph{syntactic microscope}, verifying step-by-step derivations internally without relying on external metatheories. 

Within this dual-aspect framework, the Unprovability of Consistency Thesis (UCT) is dead as a foundational barrier, since formal systems can certify the consistency of their actual derivations internally. It is time to let UCT go. Yet it served an immense cognitive purpose: by dramatizing the unprovability of $\mathsf{Con}_{\mathsf{PA}}$, it exposed the subtle divide between finite syntactic objects and non-standard models, driving a century of advances in proof-theoretic strength, reflection principles, and justification logic.

\subsection{Related Approaches and Other Theories}
Harald Grobner's preprint \cite{Grobner2026} extends the framework of selector proofs from Peano Arithmetic to all uniformly reflexive arithmetizable theories, including Zermelo--Fraenkel set theory (\(\mathrm{ZF}\)). While G\"odel's Second Incompleteness Theorem blocks internal proof of a single universal consistency sentence, primitive recursive selectors can uniformly verify local consistency schema instances by extracting the finite axiom fragments actually used in derivations. 

Complementing our approach, Santos, Sieg, and Kahle \cite{SSK24} advance Hilbert's numeral-oriented consistency program through \emph{Numeral Consistency} and \emph{Partial Finitism}, using intensional axiom numerations in $\mathsf{I\Sigma }_{\mathsf{1}}$ rather than selector proofs under standard G\"odelian formalization.

An alternative approach to salvaging parts of Hilbert's program using primitive recursive operations on derivation codes was suggested in \cite{Ign94}; however, it remained within traditional relative consistency reductions rather than proving the internal consistency of \pa.

In \cite{Kushida2021},  Kushida explores the version of a constructive proof of $\forall x\, F(x)$ stating that ``for all $x$ there is a proof of $F(x)$,'' which uses the modal provability predicate rather than parameterized selectors. Kushida provided answers to all constructively true/false questions for Friedman's closed fragment of \(\mathsf{PA}\). This work also showed the foundational limits of the implicit presentation of even constructive proofs, which paved the way to the selector proofs framework.

\subsection{ATP, AI, and the Lucas--Penrose Argument}

On the practical side, refuting UCT removes a long-standing perceived logical barrier to self-verifying Automated Theorem Proving and Artificial Intelligence. In particular, it undercuts the standard Lucas--Penrose argument, which claims that human minds possess an inherent cognitive advantage over formal machines because any consistent machine \(T\) suffers from a G\"odelian ``blind spot'' regarding its own consistency.
While a machine cannot prove its over-strengthened G\"odelian sentence \(\mathrm{Con}_{T}\) (or step non-conservatively to \(T + \mathrm{Con}_T\) without altering its axioms), \(\mathrm{Con}_{T}\) is not what genuine consistency requires. Because \pa\ can formally verify its own Hilbertian consistency via an internal selector proof, formal systems face no intrinsic logical barrier preventing self-verification of their actual derivations. Once this supposed impossibility is set aside, the challenge of safe and reliable AI shifts from metaphysical limitations to concrete engineering, alignment, and formal verification.

\section{Discussion and Common Objections}
We conclude by addressing several recurring questions from foundational discussions:
\begin{itemize}
\item \emph{Does this require first settling the philosophical definition of Hilbertian finitism?} \
No. Our theorem specifically addresses the provability of \pa-consistency \emph{by means formalizable in \pa}. While this vindicates the feasibility of Hilbert's consistency program, it relies strictly on standard G\"odelian formalization in \pa\ rather than a philosophical demarcation of finitism.
\item \emph{Does the proof rely on second-order logic or assume the standard model of \pa?} \
No. The informal quantifier over finite sequences $S$ in \eqref{eq:hilbert_orig} is handled standardly by reasoning about an arbitrary G\"odel numeral $d$ inside first-order \pa. Standard and non-standard models are discussed solely at the metalevel to explain why $\mathsf{Con}_{\mathsf{PA}}$ fails to capture consistency; neither the contentual proof nor its \pa-formalization presupposes any model.
\item \emph{Does this imply that every true $\Pi_1$-sentence becomes \pa-provable as a schema?} \
No. The requirement that the selector proof be formalizable in \pa\ acts as a strict filter. For example, for $T = \pa + \mathsf{Con}_{\mathsf{PA}}$, the sentence $\forall x\mathit{Consistent}_T(x)$ is true, yet $\mathsf{ConS}_T$ is not selector-provable in \pa\ by the Kurahashi--Sinclaire theorem.
\item \emph{Isn't a selector proof merely a disguised application of Shoenfield's restricted $\omega$-rule?}  No. Shoenfield's restricted $\omega$-rule \cite{Sho69} infers the single universal sentence $\forall x F(x)$ from the external metamathematical fact that a computable function enumerates proofs of $F(n)$ for every standard numeral $n$. That rule is non-conservative: adding it to \pa\ strictly increases proof-theoretic strength, allowing the derivation of $\mathsf{Con}_{\mathsf{PA}}$. A similar non-conservative pattern appears in paracomplete Reflective Grounded Arithmetic $(\mathrm{RGA})$, where  Bryan Ford introduces an explicit structural inference rule (adoptionRGA) that strictly extends the base system to \(\mathrm{RGA}^{+}\) in order to derive its single universal consistency sentence  \cite{Ford2026}. By contrast, working within classical \pa, a selector proof $\langle s(x), v \rangle$ does not conclude $\forall x F(x)$, nor does it rely on an unformalized external verification that $s$ succeeds. Selector proofs remain completely conservative over \pa\ (Proposition~\ref{prop:conservativity}) without inflating the theory's deductive strength or needing to derive the non-standard G\"odelian sentence $\mathsf{Con}_{\mathsf{PA}}$.
\end{itemize}

 \section{Concluding Remarks} 

Our basic findings separate into two independent claims:
\begin{enumerate}
\item \textbf{Disqualifying $\mathsf{Con}_{\mathsf{PA}}$:} 
$\mathsf{Con}_{\mathsf{PA}}$ is strictly stronger over $\mathsf{PA}$ than $\mathsf{PA}$-consistency. Thus, $\mathsf{PA} \nvdash \mathsf{Con}_{\mathsf{PA}}$ does not justify UCT.
\medskip
\item \textbf{Validating $\mathsf{ConS}_{\mathsf{PA}}$:} 
$\mathsf{ConS}_{\mathsf{PA}}$ is provably in $\mathsf{PA}$ equivalent to $\mathsf{PA}$-consistency---whether we like it or not---and admits a selector proof in \pa.
\end{enumerate}

Even if one rejects $\mathsf{ConS}_{\mathsf{PA}}$, $\mathsf{Con}_{\mathsf{PA}}$ still does not work for UCT. Without $\mathsf{ConS}_{\mathsf{PA}}$, the provability of $\mathsf{PA}$-consistency in $\mathsf{PA}$ is not settled by G2---it is left \textbf{unanswered}.

Note that despite the rehabilitation of Hilbert's consistency program, its classical cross-theoretic goals---such as proving {\sf ZF}-consistency in \pa---remain out of reach. However, another principal pillar of Hilbert's program---self-verification by proving $T$-consistency in $T$---is now a ubiquitous reality with far-reaching consequences. 

\pa\ can selector-prove many known consistency-like principles {\bf about itself}: $n$-Con, Local Reflection, the totality of the Goodstein function, the Paris--Kirby hydra battle, Beklemishev's worm principle, the Paris--Harrington Principle, etc. \cite{Gadsby2025}, but not about stronger theories. 

Here is an immediate application of our proof of \pa-consistency. 
There is a well-known verification question: suppose $\pa \proves \forall x\, F(x)$ where $F(x)$ is decidable. Why does $F(n)$ hold for every standard numeral $n = 0,1,2,\ldots$ {\em without appealing to the standard model}? Here is the answer based on the proof of \pa-consistency:
\begin{enumerate}
\item By universal instantiation, $\pa \proves F(n)$ for each numeral $n$.
\item If $F(n)$ failed for some $n$, $\Sigma_1$-completeness would give $\pa \proves \neg F(n)$, yielding a concrete \pa-derivation $k$ of $(0\!=\!1)$.
\item Our selector proof certifies \pa-consistency, ruling out any such concrete derivation $k$.
\end{enumerate}
\noindent
 Via ${\sf ConS}_{\sf PA}$, \pa\ acts as a {\bf self-certified verification engine} for $\Pi_1$ claims---no external metatheory or semantic faith in the standard model is needed anymore.

\end{document}